\documentclass{iasr}
\renewcommand\thisnumber{0}
\renewcommand\thisyear {2013}
\renewcommand\thismonth{November}
\renewcommand\thisvolume{00}
\usepackage[dvips]{graphicx}
\usepackage{amsthm}

\begin{document}

\title{On the Order of the Schur Multiplier of a Pair of Finite $p$-Groups}
\author[A. Hokmabadi et.~al.]{Azam Hokmabadi\affil{1},
       Fahimeh Mohammadzadeh\affil{1}, and Behrooz Mashayekhy\affil{2}\comma\corrauth}
 \address{\affilnum{1}\ Department of Mathematic, Faculty of Sciences, Payame Noor University, 19395-4697 Tehran, Iran. \\
           \affilnum{2}\ Department of Pure Mathematics, Center of Excellence in Analysis on Algebraic Structures, Ferdowsi University of Mashhad, P. O. Box 1159-91775, Mashhad, Iran.}
 \corraddr{Behrooz Mashayekhy, Department of Pure Mathematics, Center of Excellence in Analysis on Algebraic Structures, Ferdowsi University of Mashhad, P. O. Box 1159-91775, Mashhad, Iran.
          Email: \tt bmashf@um.ac.ir}

\newtheorem{thm}{Theorem}[section]
 \newtheorem{cor}[thm]{Corollary}
 \newtheorem{lem}[thm]{Lemma}
 \newtheorem{prop}[thm]{Proposition}
 \newtheorem{defn}[thm]{Definition}
 \newtheorem{rem}[thm]{Remark}
\newtheorem{Example}{Example}[section]

\begin{abstract}

In 1998, G. Ellis defined the Schur multiplier of a pair $(G,N)$ of groups and
mentioned that this notion is a useful tool for studying pairs of groups. In this paper, we characterize the structure of a pair of finite $p$-groups $(G,N)$ in terms of the order of the Schur multiplier of $(G,N)$ under some conditions.

\end{abstract}

\maketitle

\section{Introduction}
\label{sec1}
A definition for the Schur multiplier of a group $G$ is as the
abelian group $M(G)=(R\cap F')/[R,F]$
in which $F/R$ is a free presentation  of $G$.
In 1956, J.A. Green \cite{gr} showed that the order of the Schur
multiplier of a finite $p$-group of order $p^{n}$ is bounded by
$p^{\frac{n(n-1)}{2}}$, and hence equals to
$p^{\frac{n(n-1)}{2}-t}$, for some non-negative integer $t$.
Ya.G. Berkovich \cite{ber}, X. Zhou \cite{zu}, G. Ellis \cite{elis99}, P. Niroomand \cite{pn4,pn5} and E. Khamseh et al. \cite{kha} determined
the structure of $G$ for $0\leq t\leq 6$ by different methods. \

A pair of groups $(G,N)$ is a group $G$ with a normal subgroup $N$.
In 1998, Ellis \cite{elis98} defined the Schur multiplier of a pair $(G,N)$ of groups to be the abelian
group $M(G,N)$ appearing in the natural exact sequence
\begin{eqnarray*}
 H_3(G) &\rightarrow& H_3(G/N) \rightarrow M(G,N) \rightarrow M(G)
\rightarrow M(G/N)\ \ \ \ \ \ \  (1.1) \\
&\rightarrow& N/[N,G]\rightarrow (G)^{ab} \rightarrow (G/N)^{ab} \rightarrow 0
\end{eqnarray*}
in which $H_3(G)$ is the third homology of $G$ with integer coefficients.
He \cite{elis98} mentioned that there are three reasons for considering the Schur multiplier of a pair.\\
1) It leads to a more systematic treatment of a number of elementary results on the usual multiplier $M(G)$.
In several instances this treatment yields sharper results.\\
2) It is a useful tool for studying pairs of groups. \\
3) It can provide non-trivial information on the third integral homology of a group.

In general, there is no more clear relation between the Schur multiplier of a pair of groups and the usual Schur multiplier than the exact sequence (1.1). However, if $G$ is a semidirect product of $N$ and $K$ with $N \lhd G$, then the Schur multiplier of $(G,N)$ is more available to study. In fact, there is an isomorphism  $M(G) \cong M(G,N) \times M(K)$ 
(see \cite[page 356]{elis98}).
Recently, some results have been given about the Schur multiplier of a pair of groups $(G,N)$ when $G$ is a semidirect product of $N$ by a complement $K$ 
(see \cite{karimi, sany, chiti, sal}).

Ellis \cite{elis98} also gave an upper bound for the order of the Schur multiplier of pairs of finite groups.
It is interesting to know for which classes of pairs of groups the structure of the pair
$(G,N)$  can be completely described in terms of the order of $M(G,N)$.
In 2006, Salemkar, Moghaddam and Saeedi \cite{sal} tried to answer to this question
for a pair of finite $p$-groups and
proved the following theorem.

\begin{thm}
(\cite{sal} ). Let $(G,N)$ be a pair of
groups and $K$ be a complement of
$N$ in $G$ with $|N|=p^n$ and $|K|= p^m$. Then the following statements hold:\\
$(i)$ $|M(G,N)| \leq p^{\frac{1}{2}n(2m+n-1)}$;\\
$(ii)$ If $G$ is abelian, $N$ is elementary abelian and
$|M(G,N)| = p^{\frac{1}{2}n(2m+n-1)}$, then $G$ is elementary abelian;\\
$(iii)$ If the pair $(G,N)$ is non-capable and $|M(G,N)| = p^{\frac{1}{2}n(2m+n-1)-1}$,
then $G \cong {\bf{Z}}_{p^2}$.
\end{thm}

In this paper we extend the above theorem and characterize the structure
of the pair $(G,N)$ of finite $p$-groups in terms of the order of $M(G,N)$ in more cases.
Let $(G,N)$ be a pair of finite $p$-groups and $K$ be a complement of
$N$ in $G$. In other words, $G$ is a semidirect product of $N$ and $K$. Let $|N|=p^n$, $|K|= p^m$ and
$|M(G,N)| = p^{\frac{1}{2}n(2m+n-1)-t}$. Then we prove $t=0$ if and only if $N$ is trivial or $(G,N)$ is a pair of elementary abelian $p$-groups. Also, we determine the pair
$(G,N)$, for $t=1,2,3$. The main results of this paper are somehow generalizations of the results of \cite{ber}, \cite{zu} and \cite{elis99} to the pair of finite $p$-groups.

\section{Preliminaries}
\label{sec2}
Let $(G,N)$ be a pair of groups. We recall that a relative central extension of
the pair $(G,N)$ consists of a group homomorphism $\sigma: M \rightarrow G$,
together with an action of $G$ on $M$, such that \\
$(i)$ $ \sigma (M)=N$;\\
$(ii)$ $ \sigma (m^g)=g^{-1} \sigma (m)g$, for all $g \in G$, $m \in M$;\\
$(iii)$ $m^{\sigma(m_1)}=m_1^{-1} m m_1$, for all $m, m_1 \in M$;\\
$(iv)$ $ G$ acts trivially on $\ker \sigma$.\\

The $G$-commutator subgroup of $M$ is defined to be the subgroup
$[M,G]$ generated by the $G$-commutators $[m,g]={ m^{-1}m^{g}}$  for all
$g \in G , m \in M$, in which $m^g$ is the action of $g$ on $m$, and the $G$-center of $M$ is the central subgroup
$$Z(M,G)= \{ m \in M| m^g=m \  \ for\ all\  g \in G \}.$$
Also, the subgroup $Z_2(M,G)$ is defined by
$$\frac{Z_2(M,G)}{Z(M,G)}=Z(\frac{M}{Z(M,G)},G).$$

A pair $(G,N)$ is said to be capable if it admits a relative
central extension $\sigma: M \rightarrow G$ with $\ker \sigma=Z(M,G)$.
Note that a group $G$ is capable precisely when the pair $(G,G)$ is capable.\

We call a pair $(G,N)$ an extra special pair of $p$-groups when $Z(N,G)$ and
$[N,G]$ are the same subgroups of order $p$.
Also, we need to recall the definition of a covering pair.

\begin{defn}
(\cite{elis98} ). A relative central
extension $\sigma: N^*\rightarrow G$ of a pair $(G,N)$ will be
called a covering pair if there exists a subgroup $A$ of $N^*$
such that \\
$(i)$ $A \leq Z(N^*,G)\cap [N^*,G]$;\\
$(ii)$ $A \cong M(G,N)$;\\
$(iii)$ $N \cong N^*/A.$
\end{defn}
Note that Ellis proved that every pair of finite groups has a covering pair (see [2, Theorem 5.4]).
The following theorem plays an important role to prove the main results.

\begin{thm}
(\cite{elis98} ). Let $(G,N)$ be a pair of groups and
$K$ be a complement of $N$ in $G$. Then $$M(G) \cong M(G,N) \times M(K).$$
\end{thm}

We recall from \cite{kar} that if $G= N \times K$, then
$$M(G)\cong M(N) \times M(K) \times (N^{ab}\otimes K^{ab}).$$

A useful consequence of this fact is stated in the following.

\begin{cor}
If $G= N \times K$, then
$$|M(G,N)|=|M(N)||N^{ab}\otimes K^{ab} |.$$
\end{cor}

The following theorems will  be used in the next section. Here $D$ denotes the
dihedral group of order 8, $Q$ denotes the quaternion group of order 8 and $E_1$
and $E_2$ denote the extra special groups of order $p^3$ with odd exponent $p$ and $p^2$, respectively.

\begin{thm}
(\cite{elis99} ). Let $G$ be a group of order $p^n$.
Suppose that $|M(G)|= p^{\frac{1}{2}n(n-1)-t}$. Then \\
$(i)$ $t=0$ if and only if $G$ is elementary abelian (\cite{ber});\\
$(ii)$ $t=1$ if and only if $G \cong {\bf{Z}}_{p^2}$ or $G \cong E_1$ (\cite{ber});\\
$(iii)$ $t=2$ if and only if $G \cong {\bf{Z}}_{p} \times {\bf{Z}}_{p^2}$, $G\cong D$
or $G\cong {\bf{Z}}_{p}\times E_1$ (\cite{zu});\\
$(iv)$ $t=3$ if and only if $G\cong {\bf{Z}}_{p^3}$,
$G\cong {\bf{Z}}_{p}\times {\bf{Z}}_{p} \times {\bf{Z}}_{p^2}$, $G\cong Q$, $G\cong E_2$,
$G\cong D \times {\bf{Z}}_{2}$ or $G \cong E_1 \times {\bf{Z}}_{p} \times {\bf{Z}}_{p}$.
\end{thm}

\begin{thm}
(\cite{kar} ). Let $G \cong
{\bf{Z}}_{n_1}\times {\bf{Z}}_{n_2}\times ...\times {\bf{Z}}_{n_k}$,
where $n_{i+1}|n_i$ for all $i \in {1,2,...,k-1}$ and $k\geq 2$,
and let ${\bf{Z}}_{n}^{(m)}$ denote the direct product of $m$
copies of ${\bf{Z}}_{n}$. Then $$M(G)\cong
{\bf{Z}}_{n_2}\times {\bf{Z}}_{n_3}^{(2)}\times ...\times {\bf{Z}}_{n_k}^{(k-1)}.$$
\end{thm}

\begin{thm}
(\cite{elis98} ). Let $(G,N)$ be a pair of groups such
that $N/Z(N,G)$ is finite of prime power order $p^n$ and $G/N$ is
finite of prime power order $p^m$. Then $$ |M(G,N)||[N,G]|\leq
p^{n(2m+n-1)/2}.$$
\end{thm}

\begin{thm}
(\cite{sal} ). Let $(G,N)$ be a
pair of finite $p$-groups with $G/N$ and $N/Z(N,G)$ of orders $p^m$ and
$p^n$, respectively. If $|[N,G]|=p^{\frac{1}{2}n(2m+n-1)}$, then either
$N/Z(N,G)$ is elementary abelian or the pair
$(G/Z(N,G) ,N/Z(N,G))$ is an extra special pair of finite
$p$-groups.
\end{thm}

\begin{thm}
(\cite{sal} ). Let $(G,M)$ be a
pair of groups with $G/M$ and $M/Z(M,G)$ of orders $p^m$ and
$p^n$, respectively. Suppose $z \in  Z_2(M,G)-Z(M,G)$ and
consider two non-negative integers $\mu (z) $ and $\nu(z)$, where
$$p^{\mu (z)}=|[G,z]|\ \  , p^{\nu (z)}=|\frac{G/[G,z]}{Z(G/[G,z]),
M/[G,z])}|.$$ Then \\
$(a)$ $|[M,G]|\leq p^{\frac{1}{2}(\nu(z)(\nu(z)-1)-m(m-1))+\mu (z)} \leq
p^{\frac{1}{2}n(2m+n-1)}$.\\
$(b)$ Suppose for some non-negative integer $s$,
$|[M,G]|=p^{n(2m+n-1)/2-s}$, then the following statements hold:\\
$(i)$ $|[M/Z(M,G),G/Z(M,G)]| \leq p^{s+1}$. If $|[M/Z(M,G),G/Z(M,G)]|=
p^{s+1-k}$ for some $0 \leq k \leq s+1$, then $exp(Z_2(M,G) /
Z(M,G))\leq p^{k+1}$ and $ m+n-1-s \leq \mu (z) \leq m+n-1-s+k$.\\
$(ii)$ If  $exp(Z_2(M,G) /
Z(M,G)) \geq p^k$, then $m+n\leq s/(k-1)+k/2$.
\end{thm}

\section{Main Results }
In this section we always assume that $(G,N)$ is a pair
of finite $p$-groups such that $K$ is a complement of
$N$ in $G$, with $|N|=p^n$ and $|K|= p^m$, and hence $|M(G,N)|= p^{\frac{1}{2}n(2m+n-1)-t}$, for some $t\geq 0$.

Salemkar, Moghaddam and Saeedi \cite{sal} proved that if $t=0$ and $G$ is abelian and N is elementary abelian, then $G$ is elementary abelian. The first main result of this paper gives a suitable extended version of the above result similar to Berkovich's one \cite{ber}.

\begin{thm} With the previous assumptions and notation, $t=0$ if and only if $N$ is trivial or $(G,N)$ is a pair
of elementary abelian $p$-groups.
\end{thm}

\begin{proof}
Using Theorem 2.2, we have $|M(G,N)|= |M(G)|/|M(K)|$. Hence necessity is immediate
by Theorems 2.4.
For sufficiency, let the relative central extension
$\sigma: N^*\rightarrow G$ be a covering pair of $(G,N)$. Note that by \cite[Theorem 5.4]{elis98} the finite pair $(G,N)$ has at least one covering pair. Then
there exists a subgroup $A$ of $N^*$ such that
$ A \leq Z(N^*,G)\cap [N^*,G]$, $A \cong M(G,N)$ and $N \cong N^*/A$.
It is easy to see that for any $k \in K$, the map $\varphi_k : N^* \rightarrow N^*$,
defined by $\varphi_k(n^*)={n^*}^k$, is an automorphism of $N^*$. Therefore, using the homomorphism
$\psi:K \rightarrow Aut(N^*)$, given by $\psi(k)= \varphi_k$, we can define a semidirect
product of $N^*$ by $K$, denoted by $G^*$. It is straightforward to check that the
subgroups $[N^*,G]$ and $Z(N^*,G)$ are contained in
$[N^*,G^*]$ and $Z(N^*,G^*)$, respectively. Then the map $\delta : G^*
\rightarrow G$, given by $\delta(n^*k)=\sigma(n^*)k$, for $n^* \in N^*$, $k\in K$, is an
epimorphism with $\ker(\delta)=\ker(\sigma)$.
Therefore
$$|\frac{N^*}{Z(N^*,G^*)}|\leq
|\frac{N^*}{Z(N^*,G)}|\leq|\frac{N^*}{A}|=|N|=p^n$$ and $$
|\frac{G^*}{N^*}|=|\frac{G^*/A}{N^*/A}|=|\frac{G}{N}|=|K|=p^m.$$
So $|[N^*,G^*]|\leq p^{\frac{1}{2}n(2m+n-1)} $ by Theorem 2.6. This
implies that
$$p^{\frac{1}{2}n(2m+n-1)}=|M(G,N)|=|A| \leq |[N^*,G^*]|\leq p^{\frac{1}{2}n(2m+n-1)}.$$ Thus
$A=[N^*,G^*]$. It follows that $N \leq Z(G)$ and we have $G= N \times K$. Then Corollary 2.3 implies that $p^{\frac{1}{2}n(2m+n-1)}\leq p^{\ \frac{1}{2}n(n-1)}
 |N^{ab} \otimes K^{ab}|$. Hence $p^{mn} \leq |N^{ab} \otimes K^{ab}| \leq p^{md(N)}$, where $d(N)$ is the minimum number of generators of $N$.
Therefore $n=d(N)$ and hence
$N$ is an elementary abelian $p$-group and $|N^{ab} \otimes K^{ab}|= p^{mn}$. If $n=0$, then $N$ is the trivial subgroup. If $n>0$, then we have $p^{nm}=|{\bf {Z}}_{p}^{(n)} \otimes K^{ab}|= p^{d(K^{ab})n}$ and it follows that $d(K^{ab})=m$. Therefore $d(K)=m$ and so
$K$ is an elementary abelian $p$-group too. This completes the  proof.
\end{proof}

\begin{lem}
Let $(G,N)$ be a pair of $p$-groups such that $[N,G] \neq 1$. Then $Z(N,G) \cap [N,G] \neq 1$.
\end{lem}

\begin{proof}
Using the fact that $Z(N,G)= Z(G)\cap N$, the result follows.
\end{proof}

Salemkar, Moghaddam and Saeedi \cite{sal} proved that if $t=1$ and $(G,N)$ is non-capable,
then $G \cong {\bf{Z}}_{p^2}$. The second main result of this paper gives a vast generalization of  the above result similar to Berkovich's one \cite{ber}.

\begin{thm}
 With the previous assumptions if $t=1$, then one of the following cases holds:\\
$(i)$ $G\cong N \times K$ where $N \cong {\bf {Z}}_{p^{2}}$ and $K=1$;\\
$(ii)$ $G \cong N \times K$, where $N \cong {\bf {Z}}_{p}$ and $K$ is any group with  $d(K)=m-1$;\\
$(iii)$ $(G,N)$ is an extra special pair of groups which is capable.
\end{thm}

\begin{proof}
Choose $N^*, G^*$ and $A$ as in the proof of Theorem 3.1. We divide the proof in two cases:\\
\textbf{Case 1}. Suppose $A \neq Z(N^*,G^*)$.
Then $|N^*/ Z(N^*,G^*)|<|N^*/ A|=p^n$ and
$|N^*/Z(N^*,G^*)|\leq p^{n-1}$. Thus by Theorem 2.6
$$p^{\frac{1}{2}n(2m+n-1)-1} = |A| \leq |[N^*,G^*]| \leq p^{\frac{1}{2}(n-1)(2m+n-2)}.$$
This implies that $n+m \leq 2$. Since $t=1$, by Theorem 3.1 we have $n\neq 0$ and $G$ is not an elementary abelian $p$-group. Therefore we have $G = N \cong {\bf {Z}}_{p^{2}}$.\\
\textbf{Case 2.} Suppose $A = Z(N^*,G^*)$.
Then $|N^*/Z(G^*,N^*)|=p^n$. By Theorem 2.6, we have
$$p^{\frac{1}{2}n(2m+n-1)-1}= |A| \leq |[N^*,G^*]|=p^{\frac{1}{2}n(2m+n-1)-s}$$ for some $s \geq 0$.
It follows that $s\leq1$.

(i) First, assume that $s=0$. Hence Theorem 2.8 implies that
$|[N,G]|=p$ and $\exp(Z(N,G))=p$. If $Z(N,G)$ is cyclic, then
$[N,G]=Z(N,G)$ has order $p$ by Lemma 3.2 and hence $(G,N)$ is an
extra special pair of groups which is capable. If $Z(N,G)$ is not
cyclic, then $Z(N,G)=Z_2(N^*,G^*)/Z(N^*,G^*)$ has two distinct
subgroups of order $p$. Therefore there exist elements $y_0, z_0
\in Z_2(N^*,G^*) - Z(N^*,G^*)$ such that
$$|<y_0Z(N^*,G^*)>|=|<z_0Z(N^*,G^*)>|=p$$ and $<y_0Z(N^*,G^*)> \cap
<z_0Z(N^*,G^*)>=1$. Using the proof of Theorem 2.8, we have $|[G^*,y_0]|=|[G^*,z_0]|=p^{m+n-1}$.
On the other hand, $[G^*,y_0] \cong G/ C_{G^*}(y_0)$ and $[G^*,z_0] \cong G/ C_{G^*}(z_0)$.
So $|C_{G^*}(y_0)|=|C_{G^*}(z_0)|=p$. Thus we have
$$ Z(N^*,G^*) \leq C_{G^*}(y_0)\cap C_{G^*}(z_0)=[N^*,G^*].$$
This implies that $|[N,G]|=|[N^*/Z(G^*,N^*),G^*/Z(G^*,N^*)]|=1$
which is a contradiction.

(ii) Now assume that $s=1$, then $A=[N^*,G^*]$. This implies that
$N\leq Z(G)$. Thus $G \cong N \times K$. Then by Corollary 2.3 we
have $|M(G,N)|=|M(N)||N^{ab} \otimes K^{ab}|$. Thus
$p^{\frac{1}{2}n(2m+n-1)-1}\leq p^{ \frac{1}{2}n(n-1)} \times
|N^{ab} \otimes K^{ab}|$. Hence $p^{mn-1} \leq |N^{ab} \otimes
K^{ab}| \leq p^{md(N)}$. Therefore $m(n-d(N))\leq 1$. Then $m=0$
or $n-d(N)=0$ or $m=n-d(N)=1$. If $m=0$, then $G=N \cong {\bf
{Z}}_{p^{2}}$ by Theorem 2.4. If $d(N)=n$, then $N \cong {\bf
{Z}}_{p}^{(n)}$. Hence  $ p^{d(K)n}=|N^{ab} \otimes K^{ab}|=
p^{mn-1}$ and $n(m-d(K))=1$. Therefore $n=1$ and $d(K)=m-1$. In
other words, $N \cong {\bf {Z}}_{p}$ and $K$ is any group with
$d(K)=m-1$. Finally, if $m=n-d(N)=1$, then $K \cong {\bf
{Z}}_{p}$ and $N \cong  {\bf {Z}}_{p^2} \times {\bf
{Z}}_{p}^{(n-2)}$ which by Corollary 2.3 and Theorem 2.5 we must
have $n=1$ which is a contradiction.
\end{proof}

With an additional condition we will be able to state the reverse of the above theorem as follows.

\begin{thm}
Let $G \cong N \times K$. Then $t=1$ if and only if one of the following conditions holds:\\
$(i)$ $N \cong {\bf {Z}}_{p^{2}}$ and $K=1$;\\
$(ii)$ $N \cong {\bf {Z}}_{p}$ and $K$ is any group with  $d(K)=m-1$;\\
$(iii)$ $N \cong E_1$ and $K=1$.
\end{thm}

\begin{proof}
By the assumption, $G \cong N \times K$. Hence by Corollary 2.3 we have \begin{equation}
|M(G,N)|=|M(N)||N^{ab} \otimes K^{ab}|.
\end{equation}
Necessity is clear by the above equality. For sufficiency, first we suppose that $N$ is an elementary abelian $p$-group. Then $N \cong {\bf {Z}}_{p^n}$ and $|N^{ab} \otimes K^{ab}|=|M(G,N)|/|M(N)|=p^{mn-1}$. On the other hand, we have $|N^{ab} \otimes K^{ab}|=|{\bf {Z}}_{p}^{(n)} \otimes K^{ab}|=p^{nd(K)}$. Therefore $mn-1=nd(K)$ and so $n(m-d(K))=1$. This implies that $N\cong {\bf {Z}}_{p}$ and $d(K)=m-1$.

Now suppose that $N$ is not elementary abelian. Then using (3.1) we have $$|M(G,N)|< p^{\frac{1}{2}n(n-1)} |K|^{d(N)} \leq p^{\frac{1}{2}n(n-1)+md(N)}.$$ Hence $t=1$ implies that $m(n-d(N))=0$. But $n\neq d(N)$. Therefore $m=0$ and $|M(G,N)|=|M(N)|=p^{ \frac{1}{2}n(n-1)-1}$. Hence by Theorem 2.4 we have $K=1$ and $N\cong {\bf {Z}}_{p^{2}}$ or $N \cong E_1$. This completes the proof.
\end{proof}

The two following theorems determine the pair $(G,N)$, for $t=2$.

\begin{thm}. If $t=2$ and $K$ is not a normal
subgroup of $G$, then one of the following cases holds:\\
$(i)$ $G\cong D$ and  $N \cong {\bf {Z}}_4$ ;\\
$(ii)$ $G\cong E_2$ and $N \cong {\bf {Z}}_p$;\\
$(iii)$ $(G,N)$ is an extra special pair of groups which is
capable;\\
$(iv)$ $(G,N)$ is a pair of groups with cyclic center and
$|Z(G,N)|=p$ and $|[N,G]|=p^2$;\\
$(v)$ $(G,N)$ is a pair of groups with $exp(Z(G,N))= p^2$ and
$|[N,G]|=p$.
\end{thm}

\begin{proof} Choose $N^*, G^*$ and $A$ as in Theorem
3.1. Now
consider two cases:\\
\textbf{Case 1}. Suppose $A \neq Z(G^*,N^*)$. Then similar to
Theorem 3.3 one can show that $n+m \leq 3$. Since $K$ is not normal, we have $m,n\neq 0$.\\
If $m=1$ and $n=2$, then $|M(G)|=|M(G,N)|$. But $|M(G)|=p^{3-t'}$
for some $t'\geq 0$ and $|M(G,N)|=p$. Therefore Theorem 2.4 implies that $G=D={\bf {Z}}_4 > \!\!\!\!
\lhd {\bf {Z}}_2$.\\
If $m=2$ and $n=1$, then $K={\bf {Z}}_{p^{2}}$ or ${\bf {Z}}_{p}
\times {\bf {Z}}_{p}$. In the first case we have
$|M(G)|=|M(G,N)|$. But $|M(G)|=p^{3-t'}$ for some $t'\geq 0$.
Therefore Theorem 2.4 implies that $G=E_2$. In the second case we
have $|M(G)|=p |M(G,N)| $. Now by Theorem 2.4 we get a contradiction. The other cases do not happen.\\
\textbf{Case 2.} Suppose $A \cong Z(N^*,G^*)$. Then by Theorem
2.6, we have
$$p^{1/2n(2m+n-1)-2}= |A| \leq |[N^*,G^*]|=p^{1/2n(2m+n-1)-s},$$ for some $s \geq 0$.
It follows that $s\leq 2$.\\
Let $s=0$. Then similar to the proof of Theorem 3.3 (i), we can
show that $(G,N)$ is an extra special pair of groups which is
capable.\\
Let $s=1$. Then by Theorem 2.8, $|[N,G]| \leq p^2$ and
$exp(Z(G,N))=p$. If $|[N,G]| = p^2$ and $Z(G,N)$ is cyclic, then
$|Z(G,N)|=p$ and if $Z(G,N)$ is not cyclic, then similar to the proof
of Theorem 3.3 (i) we can show that $Z(G^*,N^*)= 1$ which is a
contradiction.\\
Now suppose that $|[N,G]|=p$. Then by Theorem 2.8,
$exp(Z(G,N))\leq p^2$. \\
If $exp(Z(G,N))= p$ and $Z(G,N)$ is cyclic,
then $(G,N)$ is an extra special pair of groups which is capable.
If $exp(Z(G,N))= p$ and $Z(G,N)$ is not cyclic, then using a similar method to
Theorem 3.3 (i), we get a contradiction.\\
Also $(G,N)$ may be pair of groups with $exp(Z(G,N))= p^2$ and
$|[N,G]|=p$.
\\
Let $s=2$. Then $A=[N^*,G^*]$. It follows that $N\leq Z(G)$ and
we have $G=N\times K$  which is a
contradiction.\\
\end{proof}

The next result is somehow a generalization of the Zhou's  one \cite{zu} to the pair of finite $p$-groups. In this theorem, we suppose that $K$ is normal.

\begin{thm}
Let $G \cong N \times K$. Then $t=2$ if and only if one of the following conditions holds:\\
$(i)$ $N \cong {\bf {Z}}_{p} \times {\bf {Z}}_{p^2}$ and $K=1$;\\
$(ii)$ $ N \cong D$ and $K=1$;\\
$(iii)$ $N \cong {\bf {Z}}_{p} \times E_1$ and $K=1$;\\
$(iv)$ $N \cong  {\bf {Z}}_{p^2}$ and $K \cong  {\bf {Z}}_{p}$;\\
$(v)$ $N\cong E_1$ and $K\cong {\bf {Z}}_{p}$;\\
$(vi)$ $N \cong  {\bf {Z}}_{p}\times {\bf {Z}}_{p}$ and $K$ is any group with $d(K)=m-1$;\\
$(vii)$ $N \cong  {\bf {Z}}_{p}$ and $K$ is any group with $d(K)=m-2$.
\end{thm}
\begin{proof}
Since $G\cong N \times K$, by Corollary 2.3 the equality (3.1) holds and the necessity follows.
For sufficiency, we proceed as in Theorem 3.4. First suppose that $N$ is an elementary abelian $p$-group. Then $|N^{ab} \otimes K^{ab}|=|M(G,N)|/|M(N)|= p^{mn-2}$.
On the other hand, we have $|N^{ab} \otimes K^{ab}|=p^{nd(K)}$.
Hence $mn-2=nd(K)$ and so $n(m-d(K))=2$. It follows that $n=1$ or $n=2$.
If $n=1$, then $N \cong {\bf {Z}}_{p}$ and $K$ is any group with $d(K)=m-2$.
If $n=2$, then $N \cong  {\bf {Z}}_{p} \times {\bf {Z}}_{p}$ and $K$ is
any group with $d(K)=m-1$.

Now suppose that $N$ is not an elementary abelian $p$-group.
Then (3.1) implies that
$|M(G,N)|< p^{\frac{1}{2}n(n-1)} |K|^{d(N)} \leq p^{\frac{1}{2}n(n-1)+md(N)}$.
It follows that $m(n-d(N)) <2$. Therefore $m=0$ or $m=1$.
If $m=0$, then by Theorem 2.4 we have $K=1$ and
$N \cong {\bf {Z}}_{p} \times {\bf {Z}}_{p^2}$ or $ N \cong D$ or
$ N \cong {\bf {Z}}_{p} \times E_1$. If $m=1$, then
$K \cong {\bf {Z}}_{p}$ and $d(N)=n-1$.
So $|N^{ab} \otimes K^{ab}|=p^{n-1}$. Therefore the equality (3.1) implies that
$|M(N)|= p^{\frac{1}{2}n(n-1)-1}$. Now using Theorem 2.4 we have $N \cong  {\bf {Z}}_{p^2}$
or $N \cong E_1$. This completes the proof.
\end{proof}

Using Theorem 2.2 of \cite{kha} and by a method similar to the proof of Theorem 3.5, we have the following theorem.

\begin{thm} If $t=3$ and $K$ is not a normal
subgroup of $G$, then one of the following cases holds:\\
$(i)$ $(G,N)$ is an extra special pair of groups which is
capable; \\
$(ii)$ $(G,N)$ is a pair of groups with $|Z(G,N)|=p,
|[N,G]|=p^2$ or $p^3$;\\
$(iii)$ $(G,N)$ is a pair of groups with
$|[N,G]|=p$ and $exp(Z(G,N))=p$ or $p^2$ or $p^3$; \\
$(iv)$ $G\cong D \times {\bf {Z}}_{2}$ and $N \cong \bf {Z}_2$ or $N \cong \bf {Z}_4 \times {\bf {Z}}_{2}$ ; \\
$(v)$ $G \cong  Q_8\times {\bf {Z}}_{2}$ and $|N|=2$ or $4$;\\
$(vi)$ $G \cong Q_{16}$ or $G \cong D_{16}$ and $N\cong \bf {Z}_2$;\\
$(vii)$ $|N|=p$ or $p^2$ and $G \cong E_2\times {\bf{Z}}_{p}, T_1, T_4, X_1, X_4, X_5$ or $X_6$ ($p\neq 2$);\\
$(viii)$ $N\cong \bf {Z}_p$ and $G \cong X_4 \times {\bf{Z}}_{p}, T_4 \times {\bf{Z}}_{p}, T_2, T_3, X_2, X_3, X_4, X_8$ or $X_9$ ($p\neq 2$);\\
where \\
$X_1= \langle a,b,c| a^{p^2}=b^p=c^p=1, [a,c]=b, [a,b]=[b,c]=1 \rangle,\\
 X_2= \langle a,b|
a^{p^2}=b^{p^2}=1, [a,b]=a^p \rangle,\\
X_3=\langle a,b| a^{p^3}=b^{p}=1, [a,b]=a^{p^2} \rangle,\\
X_4=\langle a,b,c| a^{p^2}=b^{p}=c^p=1, [b,c]=a^p \rangle,\\
X_5=\langle a,b,c| a^{9}=b^{3}=c^3=1, [a,b]=1, [a,c]=b,
c^{-1}b=a^{-3}b \rangle,\\
X_6=\langle a,b,c,d| a^{p}=b^{p}=c^p=d^p=[c,d]=b, [b,d]=a,
[a,b]=a,d]=[b,c]=[a,c]=1 (p>3) \rangle,\\
X_7=\langle a,b,c| a^{p^2}=b^{p}=c^p=1, [a,c]=b, [b,c]=1, [a,b]=a^p \rangle,\\
X_8=\langle a,b,c| a^{p^2}=b^{p}=1, c^p=a^p, [a,c]=b, [b,c]=1,
[a,b]=a^p \rangle,\\
X_9=\langle a,b,c| a^{p^2}=b^{p}=1, c^p=a^{\alpha p}, [a,c]=b,
[b,c]=1,
[a,b]=a^p (\alpha \neq 0,\neq$  and non residue mod $p$)$\rangle$,\\
$T_1=\langle a,b,c| a^{4}=b^{2}=c^2=1, [a,c]=b,  [a,b]=[b,c]=1 \rangle,\\
T_2=\langle a,b| a^{4}=b^{4}=1, [a,b]=a^2 \rangle,\\
T_3=\langle a,b| a^{8}=b^{2}=1, [a,b]=a^4 \rangle,\\
T_4=\langle a,b,c| a^{4}=b^{2}=c^2=1, [b,c]=a^2,[a,b]=[a,c]=1 \rangle.\\
$
\end{thm}

Finally, our last main result is somehow a generalization of the Ellis' one \cite{elis99} to the pair of finite $p$-groups.

\begin{thm}
Let $G \cong N \times K$. Then $t=3$ if and only if one of the following conditions holds:\\
$(i)$ $N \cong  {\bf {Z}}_{p^3}$ and $K=1$;\\
$(ii)$ $N \cong {\bf {Z}}_{p} \times {\bf {Z}}_{p} \times {\bf {Z}}_{p^2}$ and $K=1$;\\
$(iii)$ $N \cong  Q$ and $K=1$;\\
$(iv)$ $N \cong  E_2$ and $K=1$;\\
$(v)$ $N \cong  D \times {\bf {Z}}_{2}$ and $K=1$;\\
$(vi)$ $N \cong  E_1 \times {\bf {Z}}_{p}\times {\bf {Z}}_{p}$ and $K=1$;\\
$(vii)$ $N=K \cong  {\bf {Z}}_{p^2}$;\\
$(viii)$  $N \cong  {\bf {Z}}_{p^2}\times {\bf {Z}}_{p}$ and $K\cong  {\bf {Z}}_{p}$;\\
$(ix)$  $N \cong  D$ and $K\cong  {\bf {Z}}_{p}$;\\
$(x)$  $N \cong  {\bf {Z}}_{p}\times E_1$ and $K\cong  {\bf {Z}}_{p}$;\\
$(xi)$  $N \cong  {\bf {Z}}_{p^2}$ and $K\cong {\bf{Z}}_{p} \times {\bf{Z}}_{p}$;\\
$(xii)$  $N \cong  E_1$ and $K \cong {\bf{Z}}_{p} \times {\bf{Z}}_{p}$;\\
$(xiii)$ $N \cong  {\bf {Z}}_{p}$ and $K$ is any group with $d(K)=m-3$;\\
$(xiv)$ $N \cong  {\bf {Z}}_{p} \times {\bf {Z}}_{p}\times {\bf {Z}}_{p}$ and
$K$ is any group with $d(K)=m-1$.
\end{thm}

\begin{proof}
Necessity is straightforward.
The proof of sufficiency is similar to the proof
of previous theorems.
Suppose that $N$ is an elementary abelian $p$-group. Since $t=3$, we have $p^{nd(K)}=|N^{ab} \otimes K^{ab}|=p^{mn-3}.$ This implies that $n(m-d(K))=3$.
 So $n=1$ or $n=3$.
If $n=1$, then $N \cong {\bf {Z}}_{p}$ and $K$ is any group with $d(K)=m-3$.
If $n=3$, then $N \cong  {\bf {Z}}_{p}\times {\bf {Z}}_{p}\times {\bf {Z}}_{p}$
and $K$ is any group with $d(K)=m-1$.

Now suppose that $N$ is not an elementary
abelian $p$-group. Then we have
$|M(G,N)|<  p^{\frac{1}{2}n(n-1)+md(N)}$ and so  $m(n-d(N)) \leq 2$. This implies that $m=0$, $m=1$ or $m=2$.

If $m=0$, then by Theorem 2.4 we have $K=1$ and $N \cong {\bf {Z}}_{p^3}$ or
$N \cong {\bf{Z}}_{p}\times {\bf{Z}}_{p} \times {\bf{Z}}_{p^2}$ or
$N \cong Q$ or $N \cong E_2$ or $N\cong D\times {\bf{Z}}_{2}$
or $N \cong E_1 \times {\bf{Z}}_{p} \times {\bf{Z}}_{p}$.

If $m=1$, then $K \cong {\bf {Z}}_{p}$ and $d(N)=n-1$ or $d(N)=n-2$.
If $d(N)=n-1$, then $|N^{ab} \otimes K^{ab}|=p^{n-1}$ and so
$|M(N)|=|M(G,N)|/|N^{ab} \otimes K^{ab}|= p^{\frac{1}{2}n(n-1)-2}$. It follows
that $N \cong  {\bf {Z}}_{p^2}\times {\bf {Z}}_{p}$ or $N \cong  D$ or
$N \cong  {\bf {Z}}_{p}\times E_1$, by Theorem 2.4. If $d(N)=n-2$, then
$|M(N)|= p^{\frac{1}{2}n(n-1)-1}$ and hence $N \cong {\bf{Z}}_{p^2}$ or $N \cong E_1$
which is a contradiction.

If $m=2$, then $d(N)= n-1$ and $K \cong {\bf {Z}}_{p^2}$ or
$K\cong {\bf{Z}}_{p} \times {\bf{Z}}_{p}$.  If $K\cong {\bf{Z}}_{p} \times {\bf{Z}}_{p}$,
then $|N^{ab} \otimes K^{ab}|= p^{2(n-1)}$ and so $|M(N)|= p^{\frac{1}{2}n(n-1)-1}.$
Then $N \cong {\bf{Z}}_{p^2}$ or $N \cong E_1$.
Now suppose $K \cong {\bf{Z}}_{p^2}$. Then $N^{ab} \cong {\bf {Z}}_{p}^{(n-1)}$ or
$N^{ab} \cong {\bf {Z}}_{p^2}\times {\bf {Z}}_{p}^{(n-2)}$. If
$N^{ab} \cong {\bf {Z}}_{p}^{(n-1)}$, then $|N^{ab} \otimes K^{ab}|= p^{n-1}$
and so $|M(N)|= p^{\frac{1}{2}(n^2+n-4)}$ which is a contradiction.
If $N^{ab} \cong {\bf {Z}}_{p^2}\times {\bf {Z}}_{p}^{(n-2)}$, then
$|N^{ab} \otimes K^{ab}|= p^{n}$ and hence  $|M(N)|= p^{\frac{1}{2}(n^2+n-6)}$
which implies that $n \leq 2$. Therefore $N \cong {\bf {Z}}_{p^2}$.
This completes the proof.
\end{proof}

\section*{Acknowledgments}
The authors would like to thank the referee for useful suggestions to improve the present paper.\\
The first and the second authors were partially supported by a grant from Payame Noor University.


\end{document}